# Color identical pairs in 4-chromatic graphs

**Asbjørn Brændeland**

I argue that, given a 4-chromatic graph *G* and a pair of vertices {*u*, *v*} in *G*, <u>if the color of *u* equals the color of *v* in every 4-coloring of *G* then there is no planar supergraph of *G* where *u* and *v* are adjacent</u>. This is equivalent to the graph theoretical version of the Four Color Theorem, which says that every planar graph is 4-colorable. My argument focuses on the *coloring constraint* mechanisms that produce *color identical pairs* (e.g., {*u*, *v*}), and how these mechanisms affect the *adjaceability* of such pairs. Section 1 gives the relevant graph theoretical definitions, section 2 addresses *adjaceability* of vertices in planar graphs, section 3 explains *color identical pairs*, section 4 addresses *coloring constraints*, and section 5 gives the conclusive argument.

## 1. Relevant standard graph theory

A **graph** *G* is an ordered pair (*V*(*G*), *E*(*G*)) where *V*(*G*) is a set of **vertices** and *E*(*G*) is a set of **edges**. If *u* and *v* are vertices in *V*(*G*) and *uv* is and edge in *E*(*G*), we say that *u* and *v* are **adjacent**, and that *u* and *v* are the **end vertices** of *uv*.

Given two graphs *G* and *H*, if *V*(*H*) ⊆ *V*(*G*) and *E*(*H*) ⊆ *E*(*G*) then *H* is a **subgraph** of *G* and *G* is a **supergraph** of *H*.

A **path** is a graph with an ordered set of vertices ($v_1, …, v_n$) where for $k < n$, $v_k$ and $v_{k+1}$ are adjacent.

A **cycle** is a path where the last vertex is adjacent to the first.

A **face** in a graph is a region without edges, bounded by a cycle.

A **planar graph** is a graph that can be embedded in the plane by points and curves being associated with vertices and edges in such a way that the endpoints of a curve associated with an edge are associated with the end vertices of the edge, only the end points on a curve is associated with vertices, and curves only intersect at their end points.

A **coloring** of a graph *G* is the partitioning of *V*(*G*) into classes such that for every edge *xy* in *E*(*G*), *x* and *y* are not in the same class. The color set $C = \{1, …, k\}$ distinguishes the classes in a *k*-coloring. I let *red*, *green*, *blue* and *yellow* denote 1, 2, 3 and 4, respectively.

A **k-chromatic** graph is a *k*-color requiring, *k*-colorable graph.

A **k-critical** graph is a *k*-chromatic graph where every edge contributes to the graph's color demand.

The 1-, 2- and 3-critical graphs are exactly the *vertices*, the *edges* and the *odd cycles*. There is no known characterization of the 4-critical graphs.

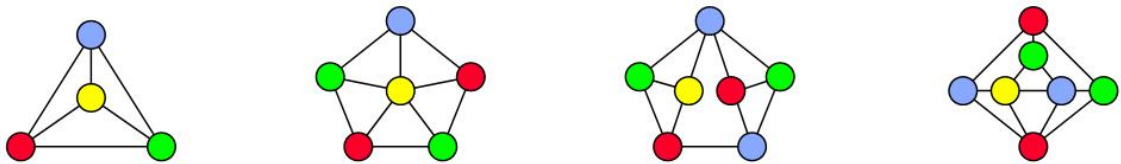

Figure 1. Four 4-critical graphs.

## 2. Adjaceability of vertices in planar graphs

**Definition 1**. Given two vertices *u* and *v* in a planar graph *G*, *u* and *v* are **adjaceable** if and only if they are not adjacent and there is no cycle *X* in *G* such that *u* and *v* are on opposite sides of *X*.

It follows that, if *uv* is not an edge in *G*, then *G* + *uv* is planar if and only if *u* and *v* are adjaceable in *G*.

The notion of *adjaceability* should be fairly intuitive, but it can also be given a theoretical basis. Let *G* be a planar graph, let *G'* be a planar embedding of *G*, let *C* be a cycle in *G* and let *C'* be the union of the curves representing the edges of *C*. Then *C'* is a closed curve, and, by *the Jordan curve theorem* (Thomassen [1992]), it divides the plane, so that if *u* and *v* are vertices in *G*, represented by *u'* and *v'* in *G'*, and *u'* and *v'* are on opposite sides of *C'*, then any curve that connects *u'* and *v'* intersects *C'*. However, since *G'* cannot contain curve intersections that do not represent vertices, there cannot be any such curve in *G'*, thus *u* and *v* cannot be adjacent in *G*.

This means that two vertices cannot be adjacent to every vertex in the same cycle on the plane without being on opposite sides of the cycle, and thus not adjaceable.

**Theorem 1**. *For every cycle X and for every vertex x on the plane, if x is adjacent to every vertex in X then there is no vertex y on the plane such that y is on the same side of X as x, and y is adjacent to every vertex in X.*

*Proof*: Let *C* be a cycle, let *a*, *b* and *c* be vertices on *C*, let *u* be a vertex adjacent to every vertex in *C*, and let *v* be a vertex on the same side of *C* as *u*. Then *v* is either separated from *a* by every cycle containing the segment *buc*, or from *b* by every cycle containing *auc*, or from *c* by every cycle containing *aub*. □

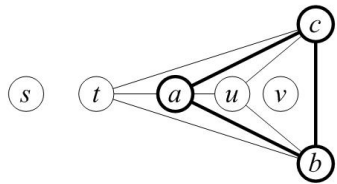

Figure 2. Here *C* = *abc* and *t* on the outside and *u* on the inside are both adjacent to all of *C*. On the outside, *s* is separated from *a* by the cycle *btc*. On the inside *v* is separated from *a* by the cycle *buc*.

**Corollary 1.** *At the most, two vertices can be adjacent to every vertex in the same cycle on the plane.*

Cycles may *cross* in the sense that one cycle can have vertices both inside and outside another. E.g. in Figure 2 the cycle *tbuc* has one vertex, *t*, outside and one vertex, *u*, inside the cycle *abc*. However, two cycles cannot cross if every vertex in both cycles is adjacent to one and the same vertex.

**Theorem 2.** *For every two cycles X and Y and for every vertex x not in V(X) ∪ V(Y), if x is adjacent to every vertex in both X and Y then X cannot have vertices both inside and outside Y.*

*Proof*: Every neighbor of *x* that is not in *V(Y)* must be on the same side of *Y* as *x*. □

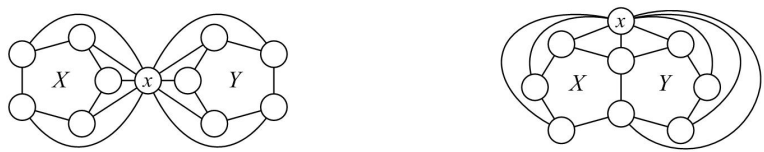

Figure 3. Since *x* is adjacent to every vertex in both *X* and *Y*, the cycles cannot *cross*, but they can *intersect*, as demonstrated by the graph to the right.

## 3. Color identical pairs

**Definition 2.** Let *u* and *v* be vertices in a *k*-chromatic graph *G*. Then {*u*, *v*} is a **color identical$_k$ pair**, or a **CI$_k$ pair**, if and only if the color of *u* equals the color of *v* in every *k*-coloring of *G*.

**Theorem 3.** *Every k-chromatic graph has a (k–1)-chromatic subgraph with an adjaceable CI$_{k-1}$ pair.*

*Proof*: Let *G* be a *k*-chromatic graph. Then there must be a subgraph *X* of *G* and an edge *xy* in *X* and a graph *Y* = *X* – *xy* such that *X* is *k*-chromatic and *Y* is (*k* – 1)-chromatic. Every *k*-critical subgraph of *G* satisfies the condition of *X*, and if there were a (*k* – 1)-coloring of *Y* in which *x* and *y* had different colors, that would also have been a coloring of *X*, and *X* would not have been *k*-chromatic. Since *x* and *y* are adjacent in *X* but not in *Y* they are adjaceable in *Y*. □

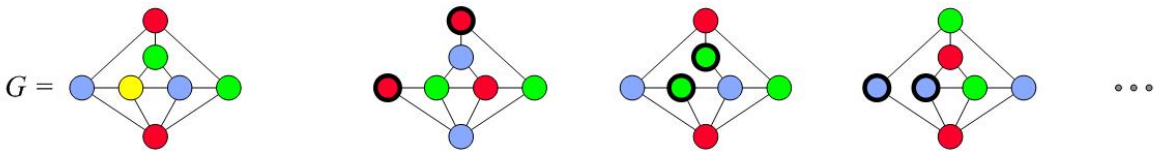

Figure 4. For every edge *xy* in the 4-critical graph *G*, (*V*(*G*), *E*(*G*) \ {*xy*}) is a 3 chromatic graph where {*x*, *y*} is a CI$_3$ pair.

The **condition of color identity** in a *k*-chromatic graph *G* must be that a pair of vertices are consistently constrained to the same color by their respective neighbors in accordance with the *coloring* definition. I.e. for every CI$_k$ pair {*x*, *y*} in *V*(*G*) there must be a pair {*X*, *Y*} of subsets of *V*(*G*), each with at least *k* – 1 vertices, such that *x* is adjacent to every vertex in *X* and *y* is adjacent to every vertex in *Y*, and for every coloring of *G* there is a set *Z* of *k* – 1 colors such that both *X* and *Y* have exactly the colors in *Z*.

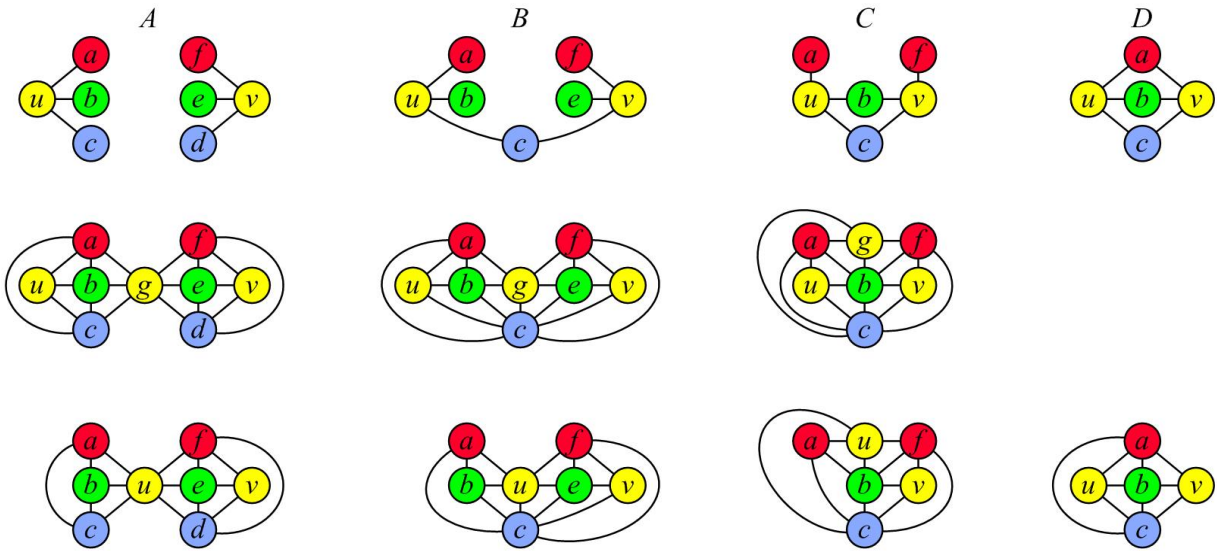

Figure 5. Four realizations of the CI$_4$ pair condition with {*u*, *v*} as the CI$_4$ pair. In *A*, {*X*, *Y*} = {{*a*, *b*, *c*}, {*d*, *e*, *f*}}, in *B*, {*X*, *Y*} = {{*a*, *b*, *c*}, {*c*, *e*, *f*}}, in *C*, {*X*, *Y*} = {{*a*, *b*, *c*}, {*b*, *c*, *f*}}, and in *D*, *X* = *Y* = {*a*, *b*, *c*}.

Figure 5 shows four realizations of the CI$_4$ condition with {*u*, *v*} as the CI$_4$ pair. The top row shows only the minimum of required edges, whereas the two rows below show examples of how the cases may be established. None of the established cases contain an adjaceable CI$_4$ pair. To determine if this lack of adjaceability is inherent to the CI$_4$ condition, we must look at the mechanisms of coloring constraints.

## 4. Coloring constraints and color fixation

### 4.1. Coloring constraint rules

Let *G* be a *k*-chromatic graph and let *R* and *S* be subgraphs of *G* such that every vertex in *S* is adjacent to every vertex in *R*. Then, in a *k*-coloring of *G*, *R constrains S* to the colors *R* does not have, and vice versa, and if the sum of the chromatic numbers of *R* and *S* equals *k*, *R fixes S* to the colors *R* does not have and vice versa.

**Definition 3**. Let *C* be the color set $\{1, \ldots, k\}$, let *G* be a *k*-chromatic graph, and let *R* be a *j*-chromatic and *S* be a $(k-j)$-chromatic subgraph of *G*. Then *S* is **color fixed** by *R* if and only if for every set of *j* colors $X \subset C$, if *R* has exactly the colors in *X* then *S* has exactly the colors in $C \setminus X$.

A color fixed graph is fixed to a *set* of colors, but for singletons I use '*v* is fixed to *c*' as short hand notation for '*v* is fixed to $\{c\}$'.

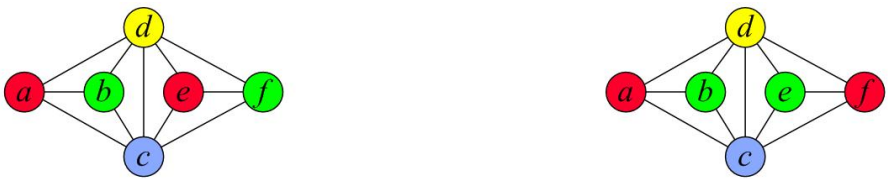

Figure 6. Two colorings of the same graph, starting with *a*, *b* and *c* (in fact the only two colorings the graph has).

Here are some of the coloring constraint relations in the graph in Figure 6:
The *red* vertex *a*
    constrains each of its neighbors *b*, *c* and *d* to {*green*, *blue*, *yellow*},
    constrains each edge *bc*, *bd* and *cd* to the same colors, and
    fixes the triangle *bcd* to the same three colors,
the *red* and *green* edge *ab*
    constrains each of *c* and *d* to {*blue*, *yellow*} and
    fixes the edge *dc* to the same colors,
the *blue* and *yellow* edge *dc*
    constrains each of *e* and *f* to {*red*, *green*} and
    fixes the edge *ef* to the same colors,
the *red*, *green* and *blue* triangle *abc*
    fixes the vertex *d* to *yellow* and
the *yellow* vertex *d*
    fixes the triangle *cef* to {*red*, *green*, *blue*}.

Constraint relations may be *indirect*. E.g., the edge *ab* fixes the edge *ef* to its own colors *indirectly*, via the edge *de*, and the triangle *abc* fixes the triangle *cef* to its own colors *indirectly*, via the vertex *d*.

### 4.2. Coloring reference

A coloring constraint relation is *mutual* but in a *description* of the coloring of a graph it can be useful to *choose* a *primary constrainer* and thus, implicitly, its *constrainees*.

**Definition 4**. Given a *k*-chromatic graph *G* we can select a (*k* – 1)-critical subgraph *R*, as ***coloring reference graph***, or ***CR-graph***, by giving *R* the first *k* – 1 colors in every *k*-coloring of *G*. *R* is then the immediate coloring constrainer of all its neighbors, the indirect coloring constrainer of every neighbor's neighbor that is not in *R* or a neighbor of *R*, and so on and so forth.

**Definition 5**. Let *G* be a graph, let *R* be a CR-graph in *G* and let *S* be a subgraph of *G* and a super-graph of *R*. Then *S* is the ***coloring constraint scope***, or ***CC-scope***, of *R* if an only if every vertex in *G* that is coloring constrained directly or indirectly by *R*, and every edge that contributes to the coloring constraint, is in *S*, and every vertex in *V*(*S*) \ *V*(*R*) is such constrained and every edge in *S* contributes to the coloring constraint.

A graph may be indirectly color fixed *as a whole* by a CR-graph *R*, without its vertices being *individually* color fixed by *R*, and a graph thus color fixed may color fix a vertex, as shown in Figure 7—where *R* is the triangle *abc*, the triangle *efh* is color fixed as a whole by *R* via *d*, and *h* is color fixed by *efg*.

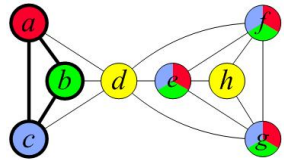

Figure 7.

### *4.3. Only an edge can color fix an edge, only an odd cycle can color fix a vertex, and only a vertex can color fix an odd cycle*

**Lemma 1.** *Given a 4-chromatic graph G, for every edge xy in G, xy is color fixed if and only if there is a color fixed edge zw in G such that both x and y are adjacent to both z and w.*

*Proof*: Let *uv* be the color fixee, let *q*, *r*, *s* and *t* be vertices such that *u* is adjacent to *q* and *r*, and *v* is adjacent to *s* and *t*, and assume that *q* and *s* are fixed to *red*, and *r* and *t* are fixed to *yellow*, as shown in Figure 8 on the top. Then each of *u* and *v* are constrained, and *uv* is fixed to {*blue*, *yellow*}.

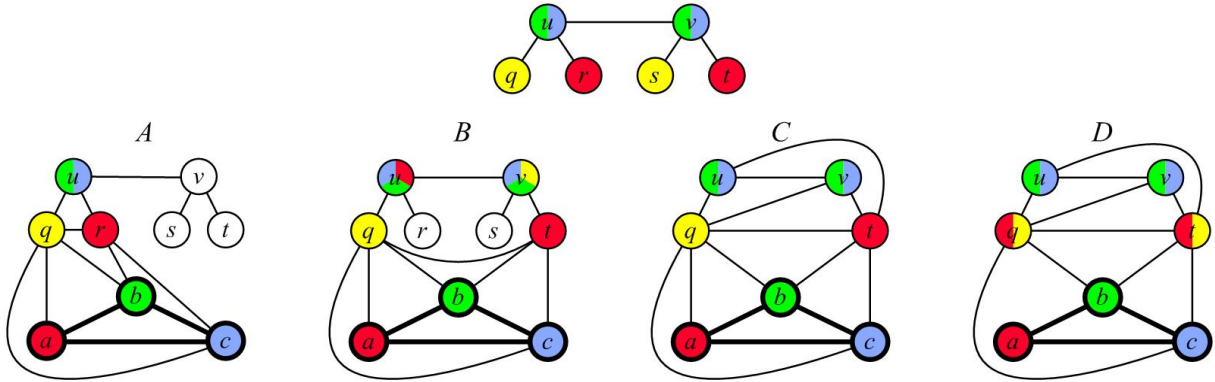

Figure 8.

However, when the triangle *abc* is introduced as CR-graph we see that it is not possible to color fix the four vertices this way. What *is* possible, is to fix one of the four vertices to *yellow* and use this and the *blue* and *green* vertices in *abc* to fix another of the four vertices to *red*, as shown in graph *C* in Figure 8, or add an edge between two of the four vertices and use the *blue* and *green* edge in *abc* to fix the added edge to {*red*, *yellow*}, as shown in graph *D* in Figure 8. In both cases we get an edge that is fixed to {*red*, *yellow*}, and this edge can in turn color fix *uv* to {*green*, *blue*}. □

**Lemma 2.** *Given a 4-chromatic graph G, for every vertex x in G, x is color fixed if and only if there is a color fixed odd cycle X in G such that x is adjacent to every vertex in X.*

*Proof*: Let *v* be the color fixee. The immediate color fixator *v* must be either (i) an *odd cycle*, or (ii) a *set of vertices*, or (iii) a set of at least one *vertex* and one *edge*.
(i) If *v* is adjacent to every vertex in a color fixed odd cycle, the cycle color fixes *v*.
(ii) If the neighbors of *v* are color fixed individually to altogether three colors, *v* is fixed to the fourth color, thus *v* and its neighbors all belong to the same *uniquely 4-colorable graph*. By Fowler [1998] a uniquely 4-colorable graph must be an *Apollonian network*, i.e., a graph that is formed by recursively subdividing a triangular face into three triangular faces, thus three of the neighbors of *v* must be the vertices of the same triangle.
(iii) Let *st* be a color fixed edge. By lemma 1 the color fixator of *st* must be another edge. Let that be *qr*. If *v* is adjacent to *s* and *t*, it must be inside one of the triangles *qst* and *rst*, and if *v* is adjacent to one of *q* and *r*, it is color fixed if *qr* is color fixed vertex by vertex, and not just as a whole. Inside the triangle there is not room for other edges or vertices that can contribute to the immediate color fixation of *v*. □

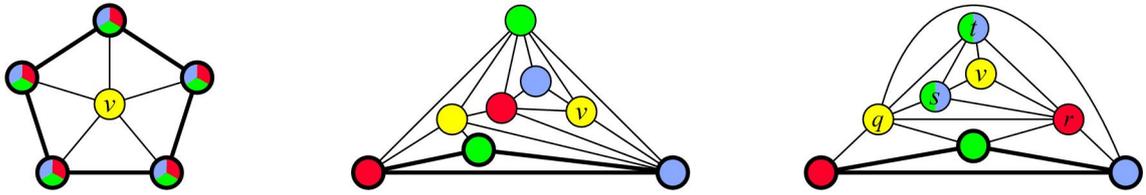

Figure 9.

**Lemma 3**. *Given a 4-chromatic graph G, for every odd cycle X in G, X is color fixed if an only if G contains a color fixed vertex that is adjacent to every vertex in X.*

*Proof*: In Figure 10, each of *A*, *B* and *C* is a subgraph of some CC-scope, and the question is how the odd cycle *S* = *abc* in each graph may be color fixed in the respective scopes. Assuming that the neighbors of *S* are already color fixed, it would seem that the three graphs cover all possible ways that *S* could have been color fixed: as a whole, by a single vertex, as in *A*; by one edge being fixed by an edge and one vertex being fixed by an odd cycle, as in *B*; or by each vertex being fixed by an odd cycle, as in *C*.

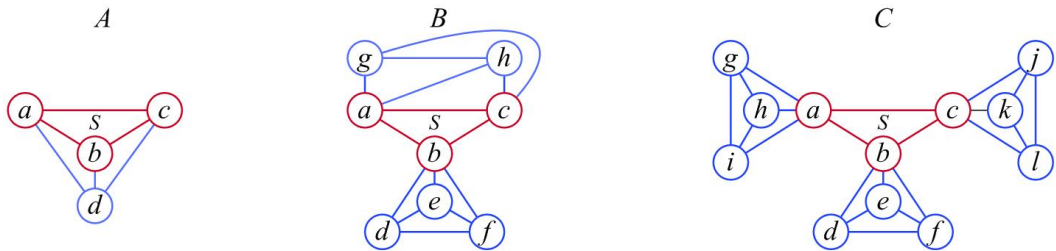

Figure 10.

In *A*, *d* may be color fixed directly or indirectly by the CR-graph that defines the CC-scope. However, in *B*, the vertices *e* and *h* are enclosed by the triangles *bdf* and *acg*, and in *C*, the vertices *e*, *h* and *k* are enclosed by the triangles *bdf*, *agi* and *cjl*, respectively, and since all the enclosing triangles intersect *S*, none of the enclosed vertices could have been color fixed without *S* already having been color fixed. The argument applies regardless of the length of *S*: A single vertex such as *d* in *A*, can be color fixed independently, and then *d* can color fix *S* as a whole, but no edge or odd cycle neighboring on *S*, such as in *B* and *C*, can be color fixed unless *S* is already color fixed.  □

## *4.4. Color fixation chains*

Color fixations *propagate* along *chains* of alternating $j$-critical and $(k - j)$-critical color fixed graphs. We are only interested in chains where $j = 1$. In a 4-chromatic graph such a chain may start with an odd cycle or a vertex and end with an odd cycle or a vertex, but with regard to color identical pairs, terminal odd cycles do not make any difference, so we will assume that the color fixation chain at hand always starts and ends with a color fixed *vertex*, as in Figure 11.

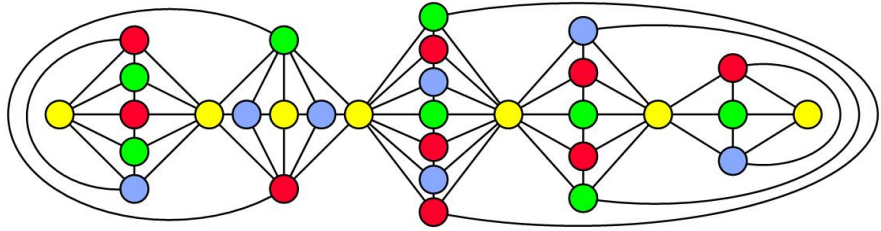

Figure 11. A 4-chromatic color fixation chain of vertices, here colored *yellow*, and odd cycles, here colored *red*, *green* and *blue*. Notice that the coloring of each of the odd cycles is arbitrary, within the given constraints.

**Definition 6**. Let $(v_1, \ldots, v_n)$ be a set of vertices and let $(C_1, \ldots, C_{n-1})$ be a set of $(k - 1)$-critical graphs. Then $(v_1, C_1, \ldots, v_i, C_i, v_{i+1}, C_{i+1}, \ldots, C_{n-1}, v_n)$ is a *color fixation$_k$ chain*, or a *CF$_k$-chain* if and only if for every $i < n$, $v_i$ and $v_{i+1}$ are directly color fixed by $C_i$, and $C_i$ and $C_{i+1}$ are directly color fixed by $v_{i+1}$.

For $k = 4$, I refer to $(v_1, \ldots, v_n)$ as the *vertex nodes* and to $(C_1, \ldots, C_{n-1})$ as the *cycle nodes* of the CF$_k$ chain.

**Lemma 4**. *Two vertices in a 4-chromatic graph are color identical if and only if they belong to the same color fixation chain.*

*Proof*: Since every pair of successive vertex nodes share an odd cycle as color fixator, and vice versa, the color fixation is repeated throughout the chain, and, since, by lemmas 2 and 3, only a vertex can color fix an odd cycle, and vice versa, color fixations cannot be transmitted between disjoint color fixation chains.  □

## 5. The non-adjaceability of CI$_4$ pairs on the plane

### *5.1. The non-adjaceability of the vertex node pairs in a CF$_4$ chain*

**Lemma 5**. *There is no pair of adjaceable vertex nodes in a 4-chromatic color fixation chain.*

*Proof*: Let $(v_1, C_1, \ldots, C_{n-1}, v_n)$ be a CF$_4$ chain. We prove the lemma by induction on $v_k$, with regard to the separation of $v_k$ from all of $(v_1, \ldots, v_{k-1})$. Assume that every vertex in $(v_1, \ldots, v_{k-2})$ are on one side, the inside, say, of $C_{k-2}$. Then, by theorem 1, $v_{k-1}$ is on the outside of $C_{k-2}$, opposite $v_{k-2}$. By theorem 2, every neighbor of $v_{k-1}$, including every vertex on $C_{k-1}$, is either *on* $C_{k-2}$, if $C_{k-2}$ and $C_{k-1}$ intersect, or *outside* $C_{k-2}$, thus everything that is inside $C_{k-2}$ is also inside $C_{k-1}$. Since $v_{k-1}$ is adjacent to both $C_{k-2}$ and $C_{k-1}$, $v_{k-1}$ lies between the two, inside $C_{k-1}$, and since $v_{k-1}$ and $v_k$ both are adjacent to all of $C_{k-1}$, again by theorem 1, $v_k$ is outside $C_{k-1}$, opposite all of $(v_1, \ldots, v_{k-1})$. We know, by theorem 1, that, since $v_1$ and $v_2$ are on opposite sides of $C_1$ they are not adjaceable, and it then follows that for every $x > 2$, $v_x$ is not adjaceable to any vertex in $(v_1, \ldots, v_{x-1})$.  □

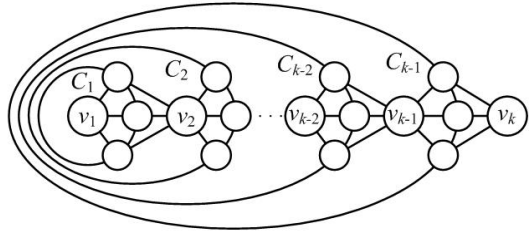

Figure 12.

**Theorem 4**. *There is no adjaceable CI$_4$ pair on the plane.*

*Proof*: By lemma 4 the members of a CI$_4$ pair must be in the same CF$_4$ chain and by lemma 5 there is no pair of adjaceable vertex nodes in a 4-chromatic color fixation chain. □

*5.2. Theorem 4 is equivalent to the Four Color Theorem*

Suppose $G$ is a planar graph with chromatic number at least 5. Without loss of generality, $G$ is 5-critical. Choose an edge $uv$ of $G$ and consider $H = G - uv$. By theorem 3 $H$ is 4-chromatic and $\{u, v\}$ is a CI$_4$ pair in $H$, but by theorem 4, $u$ and $v$ are not adjaceable on the plane, contradicting the planarity of $G$.

In the other direction, suppose there is a 4-chromatic planar graph $G$ with and adjaceable CI$_4$ pair $\{u, v\}$. Then $G + uv$ is a 5-chromatic planar graph.

**Notes**

*Color fixation chains occur wherever two or more vertices are fixed to the same color*

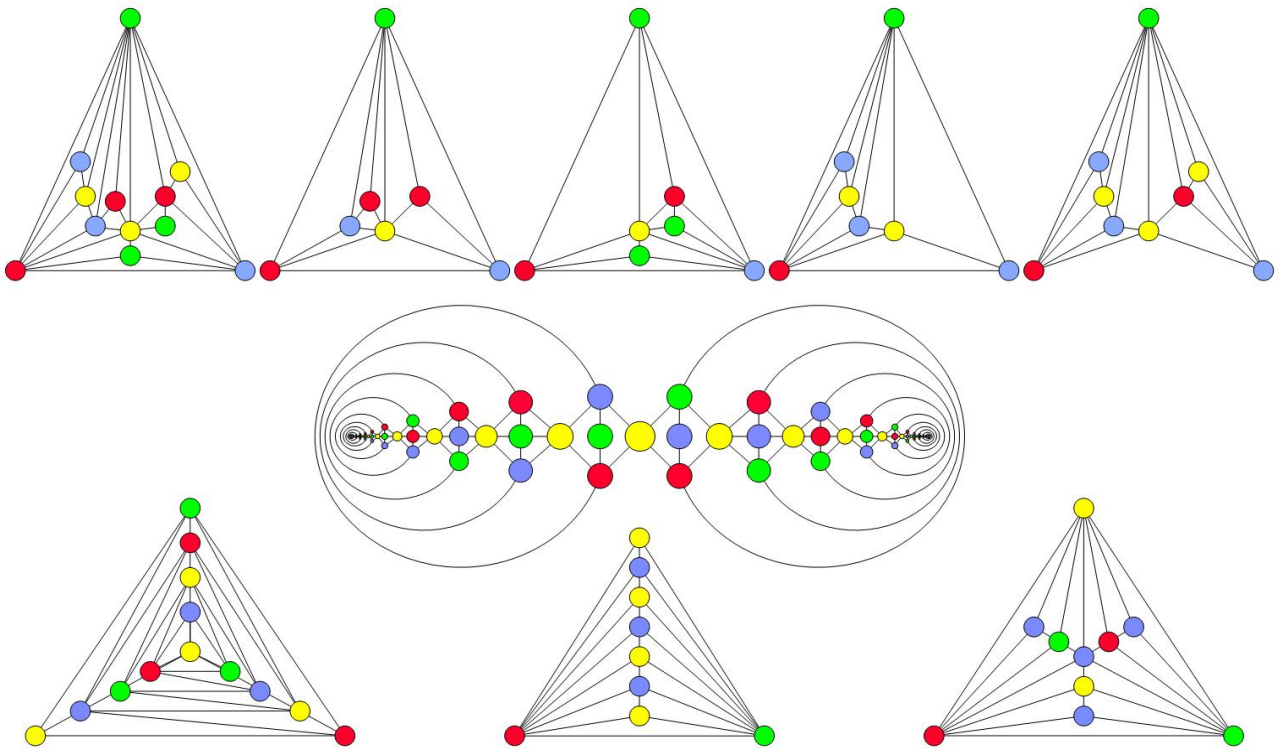

Figure 13. The top left graph contains the 4 CF$_4$ chains shown in the subsequent graphs in the top row. The middle graph is a CF$_4$ chain of its own. The bottom left graph contains 4 CF$_4$ chains. The bottom middle graph has 2 CF$_4$ chains, one with *blue* and one with *yellow* vertex nodes. The bottom right graph has 4 CF$_4$ chains, of which the one with the *blue* vertex nodes contains branches.

*Branching CF chains*

A vertex node in a $CF_4$ chain may color fix more than two cycle nodes, creating *branches*, but, since, by corollary 1, not more than two vertices can be adjacent to every vertex in the same cycle, $CF_4$ chain branching cannot occur at cycle nodes, thus the adjaceability of the vertex pairs in the chain is not affected by branching, and since, by lemma 4, every cycle node in a $CF_4$ chain separates its predecessors from its successors, a $CF_4$ chain cannot loop, so, all and all, $CF_4$ branching is harmless.

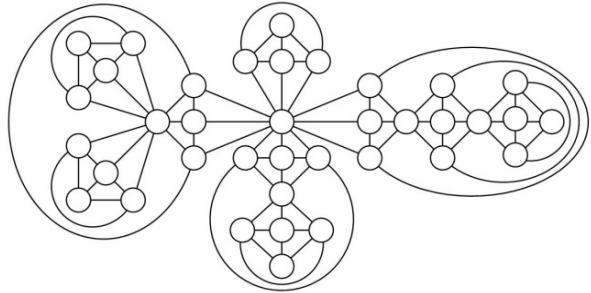

Figure 14. A $CF_4$ chain with branches.

*Inside / Outside*

The terms 'inside' and 'outside' used in the proof of lemma 4 should not be taken literally—considering the fact that the different drawings in Figure 12 and 15 all depict the same graph.

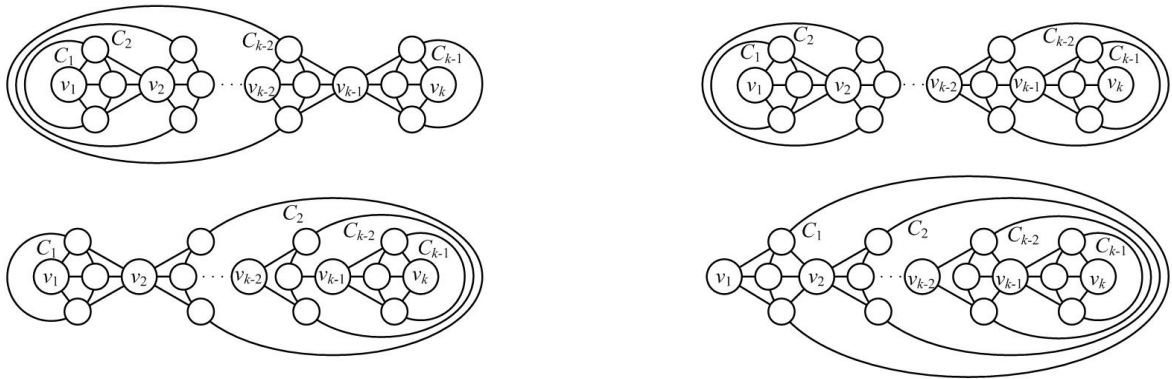

Figure 15.

*List of changes from version 2*

(not including cosmetic changes)

• The concept of *separation of vertices* has been replaced by the concept of *adjaceability of vertices*, or rather, by its negation.

• A new lemma 1 has been introduced and is used to correct an error in the proof of the previous lemma 1 (now 2).

• The previous introduction to theorem 4, explaining the equivalence of this theorem and the Four Color Theorem, has been placed in a subsection of its own, after the proof of theorem 4.